%

\magnification=1200
\input amstex
\documentstyle{amsppt}

\define\bb{{\frak b}}
\define\g{{\frak g}}
\define\m{{\Cal M([\omega]^\omega)}}
\define\bm{{B_{\Cal M}}}
\topmatter
\title
Baire Category for Monotone Sets
\endtitle
\author
Andreas Blass
\endauthor
\address
Mathematics Dept., University of Michigan, Ann Arbor, 
MI 48109, U.S.A.
\endaddress
\email
ablass\@umich.edu
\endemail
\thanks Partially supported by NSF grant DMS-9204276.
\endthanks
\abstract
We study Baire category for downward-closed subsets of $2^\omega$,
showing that it behaves better in this context than for general
subsets of $2^\omega$. We show that, in the downward-closed context,
the ideal of meager sets is prime and $\bb$-complete, while the
complementary filter is $\g$-complete.  We also discuss other cardinal
characteristics of this ideal and this filter, and we show that
analogous results for measure in place of category are not provable in
ZFC. 
\endabstract
\endtopmatter
\document

\head
1. Introduction
\endhead
We shall work with the space $2^\omega$ of infinite sequences of zeros
and ones, topologized as a product of two-point discrete spaces.
Notions of Baire category --- meager (= first category), non-meager (=
second category), and comeager (= residual) --- will always be with
respect to this topology.  We write $B$ for the $\sigma$-ideal of
meager sets in the Boolean algebra $\Cal P(2^\omega)$ of subsets of
$2^\omega$, and $B^+$ for the complement of this ideal.

We identify subsets of the set $\omega$ of natural numbers with their
characteristic functions.  Thus, we often speak of elements of
$2^\omega$ as though they were subsets of $\omega$.  We write
$[\omega]^\omega$ for the subspace of $2^\omega$ consisting of
the (characteristic functions of) infinite sets.  Restricting
attention from $2^\omega$ to $[\omega]^\omega$ does not affect Baire
category notions since the difference between the two spaces is
countable.  

We (weakly) order $2^\omega$ componentwise modulo finite sets, so
$x\leq y$ means that, for all sufficiently large $n$, $x(n)\leq y(n)$.
Under the identification of subsets of $\omega$ with their
characteristic functions, $\leq$ is the relation of
almost-inclusion ($A-B$ finite) between subsets.

Let $\m$ be the lattice of downward-closed (with respect to this
ordering) subsets $X$ of $[\omega]^\omega$.  $\m$ is ordered by
inclusion. We write $\bm$ for the $\sigma$-ideal $B\cap\m$ of meager
sets in $\m$, and we write $\bm^+$ for its complement in $\m$.  (We
use $[\omega]^\omega$ instead of $2^\omega$ only to avoid having a
non-zero intersection (namely the collection of finite subsets of
$\omega$) of all the non-zero elements of the lattice.)  The purpose
of this paper is to exhibit some pleasant properties of $\bm$ not
enjoyed by $B$.

In Section~2, we show that, for sets in $\m$, non-meagerness coincides
with groupwise density as defined in \cite5.  It follows that $\bm^+$
is a filter in $\m$ (so $\bm$ is a prime ideal) and in fact a
$(<\g)$-complete filter, where $\g$ is the groupwise density number
introduced in \cite5.

In Section~3, we study the cardinal characteristics of the ideal $\bm$
and the filter $\bm^+$.  In particular, we find that it is consistent
for the additivity of $\bm$ to be strictly greater than that of $B$.

Finally, in Section~4, we show that what we did for Baire category in
Section~2 cannot be carried over to Lebesgue measure, at least not in
ZFC alone.  Assuming the continuum hypothesis, we prove that the
complement in $\m$ of the ideal of measure-zero sets is not a filter.

\head
2. Non-Meager Sets
\endhead

This section is devoted to a combinatorial characterization of
non-meager down\-ward-closed subsets of $2^\omega$ and some immediate
consequences of this characterization.  The characterization is
essentially the same as Talagrand's characterization \cite9 of meager
filters (or ideals) on $\omega$.  Talagrand's proof does not really
use that filters are closed under finite intersections (and ideals
under finite unions).  So our proof is practically a repetition of
Talagrand's; we give it for the sake of completeness.

In \cite5, a family $G\subseteq[\omega]^\omega$ was defined to be {\sl
groupwise dense\/} if it is downward-closed (i.e., $G\in\m$) and, for
every partition of $\omega$ into finite intervals, the union of some
infinitely many of the intervals belongs to $G$.  (Actually, the
definition in \cite5 used, instead of a partition of $\omega$ into
intervals, an arbitrary family of disjoint, finite subsets of
$\omega$, but it was shown in \cite5 that no generality is lost by
considering only partitions into intervals.)

\proclaim{Theorem 1}
A set $X\in\m$ is meager if and only if it is not
groupwise dense.
\endproclaim

\demo{Proof}
Suppose first that $X$ is not groupwise dense.  Let
$\{I_0,I_1,\dots\}$ be a partition of $\omega$ into finite intervals,
no infinite union of which is in $X$.  Since $X$ is downward closed,
no element $x\in X$ includes infinitely many $I_n$.  In other words,
$X$ is covered by the sets 
$$
N_k=\{x\in\m\mid (\forall n\geq k)(\exists m\in I_n)\,x(m)=0\}.
$$
Each $N_k$ is closed and nowhere dense, so $X$ is meager.

For the converse, suppose $X\in\m$ is meager, and let it be covered by
an increasing sequence of closed, nowhere dense sets $F_n$.  We shall
partition $\omega$ into finite intervals $I_n$ and we shall define
functions $s_n:I_n\to2$ in such a way that, if $y\in2^\omega$ and the
restriction of $y$ to some $I_n$ is $s_n$, then $y\notin F_n$ for that
$n$.  Once we do this, the intervals $I_n$ will witness that $X$ is
not groupwise dense.  Indeed, suppose $x\in X$ is (the characteristic
function of) a set that includes $I_n$ for infinitely many $n$.  Then
there is $y\leq x$ whose restriction to each of these infinitely many
$I_n$ is $s_n$.  But then $y\notin F_n$ for infinitely many $n$.  As
all the $F_n$ form an increasing sequence, $y$ is not in their union
and hence not in $X$.  This is a contradiction, as $x\in X$ and $X$ is
downward-closed.

So it remains to construct the $I_n$ and $s_n$.  We proceed by
induction on $n$ and let $q$ be the smallest natural number larger
than all elements of all the previously defined $I_0$, \dots,
$I_{n-1}$.  (If $n=0$, take $q=0$.)  We shall obtain $s_n$ as the
union of $2^q$ functions $t_0,t_1,\dots,t_r$, where $r=2^q-1$ and the
domains of the $t_i$ are adjacent intervals $[q,e_0)$, $[e_0,e_1)$,
\dots, $[e_{r-1},e_r)$.  So $I_n=[q,e_r)$.  We construct $t_i$ by
induction on $i$.  

Fix a list $u_0,u_1,\dots,u_r$ of all the $2^q$ functions $[0,q)\to2$.
Then inductively choose $t_i$ so that no $y\in2^\omega$ extending
$u_i\cup t_0\cup \dots \cup t_{i-1}\cup t_i$ is in $F_n$.  Such a
choice is always possible because $F_n$ is nowhere dense and hence is
disjoint from some basic open subset of the (basic open) set of
extensions of $u_i\cup t_0\cup \dots \cup t_{i-1}$.  

Now set $s_n=t_0\cup\dots\cup t_r$.  If $y\in2^\omega$ extends $s_n$,
then, since it also extends $u_i$ for some $i$, it extends $u_i\cup
t_0\cup \dots \cup t_{i-1}\cup t_i$ and is therefore not in $F_n$.  So
$s_n$ has the required property and the proof is complete.
\qed\enddemo

The cardinal $\g$ was defined in \cite5 as the smallest number of
groupwise dense families whose intersection is empty.  This cardinal
is easily seen to be uncountable; in fact it is no smaller than the
distributivity number $\frak h$.  (We refer to \cite{10} for general
information on cardinal characteristics of the continuum and to
\cite{3, 4} for more details about $\g$.)  Thus, the intersection of any
fewer than $\g$ non-meager sets in $\m$ is non-empty.  Contrast this
with what happens in $\Cal P(2^\omega)$ where there are pairwise
disjoint, non-meager sets (in fact $2^{\aleph_0}$ of them).  The
following result slightly improves these observations by changing the
conclusion ``non-empty'' to ``non-meager.''

\proclaim{Theorem 2}
The intersection of fewer than $\g$ groupwise dense sets is groupwise
dense.  Thus, $\bm^+$ is a $(<\g)$-complete filter in $\m$.
\endproclaim

\demo{Proof}
Let fewer than $\g$ groupwise dense sets $G_i$ be given, and let $G$
be their intersection.  To show $G$ is groupwise dense, let $\omega$
be partitioned into intervals $I_n$ (listed in their natural order).
For each $i$, let $H_i$ be the collection of those infinite subsets
$x$ of $\omega$ such that $\bigcup_{n\in x}I_n\in G_i$.

We check that each $H_i$ is groupwise dense.  Let $\omega$ be
partitioned into intervals $J_k$.  Then the sets $J'_k=\bigcup_{n\in
J_k}I_n$ are intervals and constitute a partition of $\omega$.  So the
union of some infinitely many of them is in the groupwise dense family
$G_i$.  The union of the infinitely many corresponding $J_k$ is then
in $H_i$.

Since the $H_i$ are groupwise dense and there are fewer than $\g$ of
them, they have a common member $x$.  Then $\bigcup_{n\in x}I_n$ is in
all $G_i$.  So we have found a member of $G$ that is the union of
infinitely many $I_n$.
\qed\enddemo

\proclaim{Corollary}
$\bm$ is a prime ideal in $\m$.
\endproclaim

\demo{Proof}
In any lattice, an ideal whose complement is a filter is prime.
\qed\enddemo

\head
3. Cardinal Characteristics
\endhead

Associated to any ideal $I$ of sets are four cardinal characteristics,
defined as follows.  Let $X$ be the union of all the sets in $I$.

${\bold{add}}(I)$ is the smallest cardinality of a subfamily of $I$
whose union is not a member of $I$.

${\bold{cov}}(I)$ is the smallest cardinality of a subfamily of $I$
whose union is $X$.

${\bold{unif}}(I)$ is the smallest cardinality of a subset of $X$ that
is not in $I$. 

${\bold{cof}}(I)$ is the smallest cardinality of a subfamily $C$ of
$I$ such that every member of $I$ is a subset of some member of $C$.

For more information about such characteristics, see \cite{2, 6, 10}.

An ideal in $\m$, like $\bm$, is not an ideal of sets, as it is not
closed under arbitrary subsets.  Nevertheless, three of the four
cardinal characteristics make good sense in this more general context.
The exception is ${\bold{unif}}$, which would get the trivial value 1,
since singletons are not in $\bm$ simply because they are not
closed downward.  To reasonably extend ${\bold{unif}}$ to ideals in
$\m$ we define ${\bold{unif}}(I)$ to be the smallest cardinality of a
subset of $X$ not included in any member of $I$.  With this
correction, all four characteristics can also be described as the
corresponding characteristics of the ideal of sets obtained by closing
$I$ downward in $\Cal P(2^\omega)$.

Our first goal in this section is to compute the cardinal
characteristics of $\bm$.  Afterward, we shall also consider
characteristics of the complementary filter $\bm^+$.  

Recall \cite{10} that the bounding number $\frak b$ is defined to be
the smallest possible cardinality of a family $F$ of functions
$\omega\to\omega$ such that no single function $\omega\to\omega$
eventually majorizes each member of $F$.  Similarly, the dominating
number $\frak d$ is defined as the smallest cardinality of any family
$F$ of functions $\omega\to\omega$ such that every function
$\omega\to\omega$ is eventually majorized by one from $F$.

\proclaim{Theorem 3}
${\bold{add}}(\bm)={\bold{unif}}(\bm)=\frak b$, and
${\bold{cov}}(\bm)={\bold{cof}}(\bm)=\frak d$.
\endproclaim

\demo{Proof}
We record for future reference that the $X$ in the definition of
${\bold{cov}}$ and ${\bold{unif}}$, for the ideal $\bm$, is the family
of all infinite, co-infinite subsets of $\omega$.

To each partition $\Pi=\{I_0,I_1,\dots\}$ of $\omega$ into finite
intervals, associate the set $M(\Pi)$ of all infinite subsets of
$\omega$ that include only finitely many of the $I_n$.  By Theorem~1,
each $M(\Pi)$ is in $\bm$ and each member of $\bm$ is a subset of
$M(\Pi)$ for some $\Pi$.

To connect the notions involved in the definitions of the
characteristics of $\bm$ with those involved in the definitions of
$\frak b$ and $\frak d$, we use the following three constructions
relating partitions $\Pi$ as above, co-infinite subsets of $\omega$,
and functions $\omega\to\omega$.

To any $\Pi$ as above, we assign a function $F_\Pi:\omega\to\omega$ as
follows.  For any $k\in\omega$, let $n$ be the number such that $k\in
I_n$, and let $F_\Pi(k)$ be the largest element of $I_{n+2}$.

For any $f:\omega\to\omega$, let $G_f$ be an infinite, co-infinite
subset of $\omega$ such that, for each $k\in\omega$, the second
element of $\omega-G_f$ after $k$ is larger than $f(k)$.  To obtain
such a $G_f$, inductively choose the (infinitely many) elements
$a_0<a_1<\dots$ of its complement so that each $a_{n+1}$ is greater
than $f(k)$ for all $k\leq a_n$.

Finally, for any $f:\omega\to\omega$, let $H_f$ be some partition of $\omega$
into finite intervals $[a,b]$ each of which satisfies $f(a)\leq b$.
It is clear that such a partition exists; just define the intervals
one at a time by induction.

The essential properties of these constructions are given by the
following two lemmas.

\proclaim{Lemma 1}
Let $\Pi$ be a partition of $\omega$ into finite intervals, and
suppose $g:\omega\to\omega$ eventually majorizes $F_\Pi$.  Then
$M(\Pi)\subseteq M(H_g)$.
\endproclaim

\demo{Proof}
Consider any block $[a,b]$ of the partition $H_g$; by definition it
satisfies $b\geq g(a)$.  If $a$ is large enough, then by hypothesis
$g(a)\geq F_\Pi(a)$ and so $b\geq F_\Pi(a)$.  By definition of $F_\Pi$,
this means that $[a,b]$ includes an entire block of the partition
$\Pi$ (actually two entire blocks, but we don't need that here).
Therefore, a set cannot include infinitely many blocks of $H_g$
without also including infinitely many blocks of $\Pi$.  By definition
of $M(\Pi)$ and $M(H_g)$, this completes the proof.
\qed\enddemo

\proclaim{Corollary}
$\frak b\leq {\bold{add}}(\bm)$.
\endproclaim

\demo{Proof}
Let $\kappa={\bold{add}}(\bm)$.  So there are $\kappa$ sets in $\bm$
whose union is not in $\bm$.  Enlarging these sets if necessary, we
assume without loss of generality that they are $M(\Pi_i)$ for some
$\kappa$ partitions $\Pi_i$ of $\omega$ into finite intervals.  We
shall prove that $\frak b\leq\kappa$ by showing that the functions
$F_{\Pi_i}$ are not all eventually majorized by any single function
$g$.  Indeed, if $g$ eventually majorized all the $F_{\Pi_i}$, then by
the lemma the union of all the $M(\Pi_i)$ would be included in
$M(H_g)$ which is in $\bm$; this would contradict the fact that this
union is not in $\bm$.
\qed\enddemo

\proclaim{Corollary}
${\bold{cof}}(\bm)\leq\frak d$.
\endproclaim

\demo{Proof}
Let a family of $\frak d$ functions $f_i:\omega\to\omega$ be such that
every function $g:\omega\to\omega$ is eventually majorized by some
$f_i$.  Taking $g$ to be $F_\Pi$ for an arbitrary partition $\Pi$ of
$\omega$ into finite intervals, and applying the lemma, we find that
all sets of the form $M(\Pi)$ and therefore all sets in $\bm$ are
included in sets of the form $M(H_{f_i})$.  Therefore, the $\frak d$
sets of the latter form are as required in the definition of
${\bold{cof}}(\bm)$.  
\qed\enddemo

\proclaim{Lemma 2}
Let $\Pi$ be a partition of $\omega$ into finite intervals, and
suppose $f$ is such that $G_f\in M(\Pi)$.  Then $F_\Pi$ eventually
majorizes $f$.
\endproclaim

\demo{Proof}
As $G_f\in M(\Pi)$, every interval in $\Pi$, except for finitely many,
must meet the complement of $G_f$.  So, for sufficiently large $k$, as
there are two intervals of $\Pi$ between $k$ and $F_\Pi(k)$ (by
definition of $F_\Pi$), there must also be at least two elements of
$\omega-G_f$ between $k$ and $F_\Pi(k)$.  In particular, the second
element of $\omega-G_f$ after $k$ is at most $F_\Pi(k)$.  But it is
also at least $f(k)$ (by definition of $G_f$).
\qed\enddemo

\proclaim{Corollary}
$\frak d\leq{\bold{cov}}(\bm)$.
\endproclaim

\demo{Proof}
Let $\kappa={\bold{cov}}(\bm)$, and let $\kappa$ sets in $\bm$ be
given whose union contains all infinite, co-infinite subsets of
$\omega$.  Enlarging these $\kappa$ sets if necessary, we may assume
that they have the form $M(\Pi_i)$.  We shall show that the
corresponding $F_{\Pi_i}$ constitute a dominating family.  So let any
$f:\omega\to\omega$ be given.  Since $G_f$ is an infinite, co-infinite
subset of $\omega$, it lies in some $M(\Pi_i)$, and by the lemma
$F_{\Pi_i}$ eventually majorizes $f$.
\qed\enddemo

\proclaim{Corollary}
${\bold{unif}}(\bm)\leq\frak b$.
\endproclaim

\demo{Proof}
Let $\frak b$ functions $f_i$ be given, not all eventually majorized
by any single function $\omega\to\omega$.  In particular, they are not
all eventually majorized by $F_\Pi$ for any single $\Pi$.  By the
lemma, the infinite, co-infinite sets $G_{f_i}$ do not all lie in any
single $M(\Pi)$ and therefore do not all lie in any single set in
$\bm$.
\qed\enddemo

The corollaries above, together with the general facts that
${\bold{add}}\leq{\bold{unif}}$ and ${\bold{cov}}\leq{\bold{cof}}$ for
any ideal, clearly complete the proof of the theorem.
\qed\enddemo

\proclaim{Corollary}
It is consistent, relative to ZFC, that the additivity number for
$\bm$ strictly exceeds the additivity number for $B$.
\endproclaim

\demo{Proof}
In view of the theorem, this corollary merely asserts the consistency
of ${\bold{add}}(B)<\frak b$, which is well known; see for example
\cite{2, 6}.  Among the models satisfying this strict inequality are
those obtained from a model of the generalized continuum hypothesis by
adding $\aleph_2$ Laver or Mathias reals in a countable-support
iteration and the model obtained from a model of Martin's axiom and
$2^{\aleph_0}\geq\aleph_2$ by adding at least $\aleph_1$ random reals. 
\qed\enddemo

We remark that the cardinal characteristics computed for $\bm$ in
Theorem~3 are the same as the characteristics of the ideal of
$K_\sigma$ sets (countable unions of compact sets) in $\omega^\omega$.

The rest of this section is devoted to the cardinal characteristics of
the filter $\bm^+$ of non-meager sets in $\m$.  Cardinal
characteristics of a filter $F$ on a set $X$ are defined to be the
corresponding characteristics of the ideal $\{A\subseteq X|X-A\in
F\}$.  So in the case at hand we are concerned with the ideal of
non-comeager, upward-closed subsets of $[\omega]^\omega$.  If one wants
to work in the lattice $\m$ of downward-closed (rather than
upward-closed) subsets of $\omega$, one can simply replace all subsets
of $\omega$ by their complements, so the ideal becomes the ideal of
non-comeager sets in $\m$.

\proclaim{Theorem 4}
${\bold{add}}(\bm^+)={\bold{cov}}(\bm^+)=\g$ and
${\bold{unif}}(\bm^+)=2^{\aleph_0}$.
\endproclaim

\demo{Proof}
Untangling the definitions, we find that ${\bold{add}}(\bm^+)$ is the
minimum number of sets in the filter $\bm^+$ whose intersection is not
in $\bm^+$ and that ${\bold{cov}}(\bm^+)$ is the minimum number of
sets in the filter $\bm^+$ whose intersection is empty.  By Theorem~2
and the definition of $\g$, both of these cardinals equal $\g$.

${\bold{unif}}(\bm^+)$ is the minimum number of elements of $[\omega]^\omega$
needed to meet every set in $\bm^+$, so it is obviously at most
$2^{\aleph_0}$.  To prove the reverse inequality, consider any fewer than
$2^{\aleph_0}$ elements $a_i\in[\omega]^\omega$; we must find a groupwise dense
$X\in\m$ that contains none of the $a_i$.  There is an obvious choice
of $X$, namely 
$$
X=\{x\in[\omega]^\omega|(\forall i)\,a_i\not\leq x\},
$$
which contains no $a_i$ and is downward-closed.  To see that $X$ is
groupwise dense, which will complete the proof, consider an arbitrary
partition of $\omega$ into finite intervals $I_n$.  Fix a family of
$2^{\aleph_0}$ pairwise almost disjoint infinite subsets $d_\xi$ of
$\omega$.  Each $a_i$ is almost included in at most one of the sets 
$$
D_\xi=\bigcup_{n\in d_\xi}I_n,
$$ 
as the $a_i$ are infinite and the $D_\xi$, like the $d_\xi$, are
almost disjoint.  As there are more $d_\xi$'s then $a_i$'s, there must
be a $D_\xi$ that includes no $a_i$ and is therefore in $X$.  As
$D_\xi$ is a union of infinitely many $I_n$, the proof that $X$ is
groupwise dense is complete.
\qed\enddemo

We do not know the value of ${\bold{cof}}(\bm^+)$, but we have the
following partial information.

\proclaim{Theorem 5}
$\bold{cof}(\bm^+)\geq2^{\aleph_0}$ and $\bold{cof}(\bm^+)>\bb$.
\endproclaim

\demo{Proof}
The first inequality follows from Theorem~4, since
$\bold{cof}\geq\bold{unif}$ for any proper ideal.  The second
inequality follows from thie first if $\bb<2^{\aleph_0}$, so we assume
from now on that $\bb=2^{\aleph_0}$.  To show that
$\bold{cof}(\bm^+)>2^{\aleph_0}$, let $2^{\aleph_0}$ groupwise dense
families $X_\alpha$ ($\alpha<2^{\aleph_0}$) be given; we shall
construct a groupwise dense $Y$ such that no $X_\alpha\subseteq Y$.

List all the partitions of $\omega$ into finite intervals as
$\Pi_\alpha$ ($\alpha<2^{\aleph_0}$).  We construct $Y$ by an
induction of length $2^{\aleph_0}$; at each step we declare one set
$y_\alpha\in[\omega]^\omega$ to be in $Y$ and one set
$x_\alpha\in[\omega]^\omega$ to be out of $Y$.  Since $Y$ is to be
monotone, we ensure that no $y_\alpha$ is $\geq$ any $x_\beta$.  

At stage $\alpha$, we proceed as follows.  As in the proof of
Theorem~4, form $2^{\aleph_0}$ almost disjoint sets $y$, each of which
is a union of infinitely many intervals from $\Pi_\alpha$.  Each of
the $x_\beta$'s defined at earlier stages, being infinite, is $\leq$
at most one of these $y$'s.  As there are fewer than $2^{\aleph_0}$
such $x_\beta$'s, we can choose one of our $y$'s that is $\geq$ none
of them; take this $y$ as $y_\alpha$.

Next, notice that, since $\{y_\beta\mid\beta\leq\alpha\}$ has
cardinality $<2^{\aleph_0}=\bb=\bold{unif}(\bm)$, its downward closure
cannot be groupwise dense and therefore cannot include $X_\alpha$.  So
we can define $x_\alpha$ to be some element of $X_\alpha$ that is not
$\leq y_\beta$ for any $\beta\leq\alpha$.  This completes stage
$\alpha$ of our construction.

After all $2^{\aleph_0}$ stages, let $Y$ be the downward closure of
$\{y_\alpha\mid\alpha<2^{\aleph_0}\}$.  $Y$ is groupwise dense because
it contains, for each $\Pi_\alpha$, an infinite union $y_\alpha$ of
its intervals.  $Y$ includes no $X_\alpha$ since $x_\alpha\in
X_\alpha$ and $x_\alpha\notin Y$.
\qed\enddemo

\head
4. Variants 
\endhead

The two concepts connected by Theorem~1, Baire category and groupwise
density, have close relatives, to which one might reasonably try to
extend Theorem~1.  Groupwise density is, as its name suggests, a
variant of the more familiar notion of density, and Baire category is
in many respects similar to Lebesgue measure \cite7.  In this section,
we show that neither of the variants of Theorem~1 suggested by these
analogies is provable.  One fails outright by a trivial argument.  The
other is at least consistently false; we do not know whether it is
consistently true.

We treat first the easier situation, the one involving density.
Recall that a set $X\in\m$ is {\sl dense\/} if every infinite
$A\subseteq\omega$ has an infinite subset $B\subseteq A$ with $B\in
X$.  (This is the usual notion of density for the notion of forcing
consisting of the infinite subsets of $\omega$ ordered by inclusion.)
It is easy to see that groupwise density implies density; just
consider a partition of $\omega$ into intervals, each of which
contains at least one element of $A$.  It is also easy to see that the
dense sets form a filter in $\m$ and in fact a countably complete
filter; its additivity number is the cardinal $\frak h$ introduced and
studied in \cite1 and usually called the distributivity number
(because it measures the distributivity of the complete Boolean
algebra associated to the forcing mentioned above).  

A density analog of Theorem~1 would say that the filter of dense sets
is prime in $\m$, i.e., that the non-dense sets constitute an ideal
(preferably even a $\sigma$-ideal).  It is easy to see, however, that
this analog is false.  Let $X_0\in\m$ consist of those
$a\subseteq\omega$ in which all but finitely many elements are even,
and let $X_1$ be defined similarly with ``odd'' in place of ``even.''
Then neither $X_0$ nor $X_1$ is dense, but their union is dense, so
the non-dense sets fail to form an ideal.

Turning to the less trivial case of Lebesgue measure, we note that the
sets in $\m$ of measure zero constitute a $\sigma$-ideal and we ask
whether this ideal is prime, i.e., whether the sets of positive outer
measure constitute a filter (preferably a countably complete filter)
in $\m$.  The following theorem gives a consistent negative answer.
In the statement and proof of the theorem, ``measure'' refers to the
version of Lebesgue measure appropriate for the space $2^\omega$,
namely the product measure obtained from the uniform measure on 2.

\proclaim{Theorem 6}
Assume the continuum hypothesis.  Then there exist two sets in $\m$,
each of positive outer measure, whose intersection has measure zero.
\endproclaim

\demo{Proof}
Since the continuum hypothesis is assumed, let all the Borel sets of
measure 1 be listed in a sequence indexed by the countable ordinals.
We define two sequences of elements $x_\alpha$ and $y_\alpha$ of
$2^\omega$, each indexed by the countable ordinals $\alpha$, subject
to the following four requirements:
\roster
\item Both $x_\alpha$ and $y_\alpha$ belong to the $\alpha$th
Borel set of measure 1 (in our fixed list).
\item $x_\alpha$ has density 1/2 in $y_\beta$ for all $\beta<\alpha$.
\item $y_\alpha$ has density 1/2 in $x_\beta$ for all
$\beta\leq\alpha$.
\item Both $x_\alpha$ and $y_\alpha$ have density 1/2 in $\omega$.
\endroster
By ``$a$ has density 1/2 in $b$,'' where $a$ and $b$ are infinite
subsets of $\omega$, we mean that the ratio of elements of $a$ among
the first $n$ elements of $b$ tends to 1/2 as $n$ increases,
$$
\lim_{n\to\infty}\frac{|a\cap b\cap\{0,1,\dots,n-1\}|}
{|b\cap\{0,1,\dots,n-1\}|}=\frac12.
$$
We note that, by the strong law of large numbers, for any fixed
infinite $b$, almost all $a$ have density 1/2 in $b$.  If we attempt
to define the $x_\alpha$ and $y_\alpha$ by induction on $\alpha$
(defining the $x$ before the $y$ at each stage because of the $<$ in
\therosteritem2 and the $\leq$ in \therosteritem3), we find that, at
each step, the element of $2^\omega$ that we wish to define is subject
to countably many requirements, each of which is satisfied by almost
all elements of $2^\omega$.  Since measure is countably additive, the
element we need always exists (in fact, almost any element will do),
so the inductive definition succeeds.

Having defined the $x_\alpha$ and $y_\alpha$, we let $X$ be the
smallest element of $\m$ containing all the $x_\alpha$.  So $X$
consists of those $a\in[\omega]^\omega$ that are $\leq x_\alpha$ for
some $\alpha$.  Thanks to \therosteritem1, $X$ intersects every Borel
set of measure 1 and therefore has positive outer measure.  Similarly,
let $Y$ be the smallest set in $\m$ containing all the $y_\alpha$; it
too has positive outer measure.

To complete the proof, we check that $X\cap Y$ has measure zero.
Consider an arbitrary element $a\in X\cap Y$.  By definition of $X$
and $Y$, we have $a\leq x_\alpha$ and $a\leq y_\beta$ for some
$\alpha$ and $\beta$.  Consider first the case where $\alpha>\beta$.
Since $y_\beta$ has density 1/2 in $\omega$ by \therosteritem4 and
$x_\alpha$ has density 1/2 in $y_\beta$ by \therosteritem2, it follows
that $x_\alpha\cap y_\beta$ has density 1/4 in $\omega$.  The same
conclusion follows in the other case, where $\alpha\leq\beta$, by the
same argument with $x_\alpha$ and $y_\beta$ interchanged and with
\therosteritem3 in place of \therosteritem2.  In either case, $a$,
being included in $x_\alpha\cap y_\beta$ modulo a finite set, cannot
have density 1/2 in $\omega$.  This shows that $X\cap Y$ is disjoint
from the set of elements of $2^\omega$ of density 1/2.  Since the
latter set has measure 1, $X\cap Y$ has measure 0.
\qed\enddemo 

The assumption of the continuum hypothesis in Theorem~6 can be weakened
to the assumption $\bold{cov}(L)=\bold{cof}(L)$, where $L$ is the
ideal of sets of measure zero.  The only changes needed in the proof are
that the countable ordinals are replaced by the ordinals below
$\bold{cof}(L)$ and that instead of enumerating all the Borel sets of
measure 1 we enumerate only enough of them to have all the others as
supersets.

We conclude this paper with comments on some related work of Plewik
\cite8, dealing with ideals, rather than arbitrary downward-closed
families, of subsets of $\omega$.  By an ideal, we mean an $X\in\m$
such that $x\cup y\in X$ for all $x,y\in X$.  (Here we view $x$ and
$y$ as subsets of $\omega$.)  Let $\g'$ be the smallest number of
non-meager ideals with empty intersection.  Plewik \cite8 showed that
this definition is unchanged if we replace ``empty'' with meager, and
he proved $\frak h\leq\g'\leq\frak d$.  The following result describes
the connection between $\g'$ and $\g$.  It uses the {\sl splitting
number\/} $\frak s$ defined (cf. \cite{10}) as the smallest possible
cardinality for a family of subsets of $\omega$ such that every
infinite subset $x$ of $\omega$ is split by some $y$ in the family, in
the sense that both $x\cap y$ and $x-y$ are infinite.

\proclaim{Theorem 7}
$\min(\g',\frak s)\leq\g\leq\g'$.
\endproclaim

\demo{Proof}
That $\g\leq\g'$ is clear, since ideals are among the members of $\m$.
To prove the other inequality, let $\kappa<\min(\g',\frak s)$, and let
$\kappa$ groupwise dense families $G_i$ be given.  We must find an
element in their intersection.  For each $i$, let $H_i$ be the ideal
generated by $G_i$.  Being supersets of the $G_i$, the $H_i$ are
non-meager, and, since there are fewer than $\g'$ of them, all the
$H_i$ have a common member $x$.  Thus, for each $i$, some finitely
many elements $y_{ik}$ of $G_i$ cover $x$.  (Here $k$ ranges from 1 to
some finite $n_i$.)  Since the total number of all the $y_{ik}$, as
both $i$ and $k$ vary, is at most $\kappa<\frak s$, there must be an
infinite $z\subseteq x$ not split by any $y_{ik}$.  That is, for each
$i$ and $k$, either $z$ is almost included in $y_{ik}$ or they are
almost disjoint.  For any fixed $i$, the $y_{ik}$ (as $k$ varies)
cover $x$, so they cover $z$, so (because there are only finitely many
of them) they cannot all be almost disjoint from $z$.  So for each $i$
some $y_{ik}$ almost includes $z$.  But $G_i$ is downward-closed and
contains $y_{ik}$.  So $z\in G_i$ for all $i$.
\qed\enddemo

Note that the second inequality in the theorem and the fact that
$\frak h\leq\g$, pointed out in \cite3, imply Plewik's result that
$\frak h\leq\g'$.


\Refs
\ref\no 1
\by B. Balcar, J. Pelant, and P. Simon
\paper The space of ultrafilters on $N$ covered by nowhere dense sets 
\jour Fund. Math.
\vol 110
\yr 1980
\pages 11--24
\endref

\ref\no 2
\by T. Bartoszy\'nski, H. Judah, and S. Shelah
\paper The Cicho\'n diagram
\jour J. Symbolic Logic
\vol 58
\yr 1993
\pages 401--423
\endref

\ref\no 3
\by A. Blass
\paper Applications of superperfect forcing and its relatives
\inbook Set Theory and its Applications
\eds J. Stepr\B ans and S. Watson
\bookinfo Lecture Notes in Math. 1401 
\publ Springer-Verlag
\yr 1989
\pages 18--40
\endref

\ref\no 4
\by A. Blass
\paper Groupwise density and related cardinals
\jour Arch. Math. Logic
\vol 30
\yr 1990
\pages 1--11
\endref

\ref\no 5
\by A. Blass and C. Laflamme
\paper Consistency results about filters and the number of
inequivalent growth types
\jour J. Symbolic Logic
\vol 54
\yr 1989
\pages 50--56
\endref

\ref\no 6
\by D. Fremlin
\paper Cicho\'n's diagram
\inbook S\'eminaire Initiation \`a l'Analyse
\eds G. Choquet, M. Rogalski, and J. Saint-Raymond 
\publ Univ. Pierre et Marie Curie
\yr 1983/84
\pages (5-01)--(5-13)
\endref

\ref\no 7
\by J. Oxtoby
\book Measure and Category
\publ Springer-Verlag
\yr 1971
\endref

\ref\no 8
\by S. Plewik
\paper Ideals of the second category
\jour Fund. Math.
\vol 138
\yr 1991
\pages 23--26
\endref

\ref\no 9
\by M. Talagrand
\paper Compacts de fonctions mesurables et filtres non mesurables
\jour Studia Math.
\vol 67
\yr 1980
\pages 13--43
\endref

\ref\no 10
\by J. Vaughan
\paper Small uncountable cardinals and topology
\inbook Open Problems in Topology
\eds J. van Mill and G. Reed
\publ North-Holland
\yr 1990
\pages 195--218
\endref

\endRefs
\enddocument